# Some Investigations about the Properties of Maximum Likelihood Estimations Based on Lower Record Values for a Sub-Family of the Exponential Family

**Saman Hosseini[1], Parviz Nasiri[2], Dler Hussein Kadir[3], and Sharad Damodar Gore[4]**

## Abstract:

Here, in this paper it has been considered a sub family of exponential family. Maximum likelihood estimations (MLE) for the parameter of this family, probability density function, and cumulative density function based on a sample and based on lower record values have been obtained. It has been considered Mean Square Error (MSE) as a criterion for determining which is better in different situations. Additionally, it has been proved some theories about the relations between MLE based on lower record values and based on a random sample. Also, some interesting asymptotically properties for these estimations have been shown during some theories.

**Keywords**: Exponential Family, Lower Record Values, ML Estimation, Asymptotically Un-biasedness.

## 1- Introduction:

For the first time, Chandler in 1952 introduced the concept of record values and record times. After him and in a short period of time there were done a large number of researches about this new topic by statisticians. Shorrock in 1973 obtained some asymptotically conclusion in the field of record values, in the same time Ahsanullah during a paper named "an introduction to record values" clarified basic concepts of record theory (Ahsanullah 1973). Record theory is useful topic that is used in any branch of science from economy to any competitive phenomena. For the recent researches in this field it can be mentioned works of Ahmadi & Doostparast (2006), Eric Marchand & Parsian (2009), Nasiri and Hosseini (2013).

Before beginning the main topics, it is needed to introduce a sub-family of exponential family. Generally, $f(x;\theta)$ determines exponential family if

$$f(x;\theta) = h(x)g(\theta)\exp\{-\eta(\theta)T(x)\}.$$

Here, it is considered an important sub family of exponential family as bellow:
$F(x;\theta) = \exp\{-B(\theta)A(x)\}$.

[1] Correspondent author, he is a lecturer at department of computer science, Cihan University-Erbil, Kurdistan Region, Iraq. Email them at Saman.hosseini@cihanuniversity.ac.iq and s.hosseini.stat@gmail.com
[2] Parviz Nasiri is an associate professor at Payam Noor University of Tehran, Iran
[3] Dler Hussein Kadir is a lecturer at department of statistics, Salahaddin University-Erbil, Kurdistan Region, Iraq
[4] Sharad Damodar Gore is a professor at department of statistics, Savirtibai Phule Pune University, India.

It is called second type of exponential family during this paper. In which $A(x)$ is a decreasing differentiable function and

$B(\theta) > 0; x \in [a,b]; A(a) = +\infty; A(b) = 0; a,b \in \mathbb{R}$

*Table 1-1: Examples of the second type exponential family*

| $A(x)$ | $B(\theta)$ | Distribution name |
|---|---|---|
| $A(x) = x^{\alpha}$ | $B(\theta) = \theta^{\alpha}$ | Ferchet |
| $A(x) = \exp(-x)$ | $B(\theta) = \exp(\theta)$ | Gumbel |
| $A(x) = -Ln(x)$ | $B(\theta) = -Ln(x)$ | Power function |

In this section it is tried to study theorems about the relation between ML estimations based on a random sample and ML estimation based on lower record values.

## 2-Maximum Likelihood Estimations Based on Lower Record Values and Based on a Random Sample for Unknown Parameter of Type II Exponential Family

In this section it is tried to study some theorems about the relation between ML estimations based on a random sample and ML estimation based on lower record values. Also, considering MSE as medium of assessment, it is shown that in some cases estimations based on lower record values have less value of error. Remember the general form for this family is as follows.

$F(x;\theta) = \exp\{-A(x)B(\theta)\}; B(\theta) > 0; x \in [a,b]; A(a) = +\infty; A(b) = 0; a,b \in \mathbb{R}$ .

$F(x;\theta) = \exp\{-B(\theta)A(x)\}$ that’s why, probability density function is obtained as:

$f(x;\theta) = -B(\theta)A^{'}(x)\exp\{-B(\theta)A(x)\}$ .

Considering a random sample of size n, the joint probability density function for a random sample is obtained as following.

$$f(x_1,...,x_n;\theta) = \prod_{i=1}^{n} f(x_i;\theta) = \prod_{i=1}^{n} (-1)^n B^n(\theta)A^{'}(x_i)\exp\{-B(\theta)\sum_{i=1}^{n} A(x_i)\} =$$

$$B^n(\theta)\prod_{i=1}^{n}\{-A^{'}(x_i)\}\exp\{-B(\theta)\sum_{i=1}^{n} A(x_i)\}.$$

Hence,

$$L = \log f\left(x_1,\ldots,x_n;\theta\right) = nlogB\left(\theta\right) + \sum_{i=1}^{n} \log\{\frac{1}{A'(X_i)}\} - B\left(\theta\right) \sum_{i=1}^{n} A\left(X_i\right).$$

As previously said, the maximum likelihood estimation is calculated by solving the equation $\frac{dL}{d\theta} = 0$, therefore

$$\frac{\partial L}{\partial \theta} = 0 \Rightarrow \frac{nB'(\theta)}{B(\theta)} - B'(\theta) \sum_{i=1}^{n} A\left(X_i\right) = 0 \Rightarrow B\left(\theta\right) = \frac{n}{\sum_{i=1}^{n} A\left(X_i\right)}, \tag{1}$$

assuming that $B(\theta)$ is a one-to-one function, the above equation will be solved like this:

$$\theta = B^{-1}\left(\frac{n}{\sum_{i=1}^{n} A\left(X_i\right)}\right)$$

Since the cumulative distribution function $F(x;\theta) = \exp\{-A(x)B(\theta)\}$ with the properties and the conditions described for A and B is a continuous function, it is possible to derive the first and the second order derivatives. So,

$$\frac{\partial^2 L}{\partial \theta^2} = \frac{nB''(\theta)B(\theta) - n[B'(\theta)]^2}{B^2(\theta)} - B''(\theta)t, (t = \sum_{i=1}^{n} A(x_i)) \tag{2}$$

In addition, substituting (1) in (2), the below result is reached.

$$\left.\frac{\partial^2 L}{\partial \theta^2}\right|_{B\left(\theta\right) = \frac{n}{\sum_{i=1}^{n} A\left(X_i\right)}} = \frac{nB''(\theta)(\frac{n}{t}) - n[B'(\theta)]^2}{(\frac{n}{t})^2} - B''(\theta)t = -\frac{[B'(\theta)t]^2}{n} < 0 \tag{3}$$

The value of second derivative is negative therefore; the relations (1) and (3) tell us the following estimation is ML estimation.

$$\hat{\theta}_{MLE} = B^{-1}(\frac{n}{T}), T = \sum_{i=1}^{n} A\left(X_i\right)$$

On the other hand, $A(X_i)$ is distributed as Exponential distribution with parameter $\frac{1}{B(\theta)}$.

Because

$$P(A(X_i) \le t) = 1 - P(A(X_i) > t) = 1 - (1 - F(t;\theta)) = F(t;\theta) = 1 - \exp(-B(\theta)t).$$

On the other hand, $A(X_i)$'s are independent variables that's why it is possible to find the distribution of $\sum_{i=1}^{n} A(X_i)$ easily like follow:

$$T = \sum_{i=1}^{n} A(X_i) \overset{distribution}{=} Gamma(n, \frac{1}{B(\theta)}) .$$

As a result

$$\hat{\theta} \overset{distribution}{=} B^{-1}(\frac{n}{T}) \tag{4}$$

in which

$$T \overset{distribution}{=} Gamma(n, \frac{1}{B(\theta)}) .$$

On the other hand, with respect to the definition, the lower record values of the joint probability density function $(R_1^{'},...,R_n^{'})$ are resulted as follows (Arnold, Balakrishnan and Nagaraja, 1998).

$$f\left(r_1^{'},\dots,r_n^{'};\theta\right) = f\left(r_n^{'};\theta\right)\prod_{i=1}^{n-1}\frac{f\left(r_i^{'};\theta\right)}{F\left(r_i^{'};\theta\right)}$$

$$= (-A^{'}(r_n^{'})B(\theta))\exp\{-B(\theta)A(r_n^{'})\}\{\prod_{i=1}^{n-1} -A^{'}(r_i^{'})B(\theta)\}$$

$$= B^{n}(\theta)\{-\prod_{i=1}^{n} A^{'}(r_i^{'})\}\exp\{-B(\theta)A(r_n^{'})\}$$

Considering the above probability density function, the maximum likelihood function based on lower record values is obtained as:

$$L_{R^{'}} = \log(f(r_1^{'},...,r_n^{'})) = nLog\{B(\theta)\} + \sum_{i=1}^{n} Log\{\frac{1}{A^{'}(r_i^{'})}\} - A(r_n^{'})B(\theta) ,$$

the maximum likelihood estimation based on lower record values is obtained by solving the equation

$\frac{\partial L_{R^{'}}}{\partial \theta} = 0$. Therefore,

$$\hat{\theta}_{MLE-BASED-ON-LOWER-RECORDS} = \hat{\theta}_{MLE,(R_1^{'},...,R_n^{'})}$$
$$= B^{-1}(\frac{n}{A(R_n^{'})}),$$

On the other hand, $R_n^{'}$ is distributed as the following probability density function (Arnold, Balakrishnan and Nagaraja, 1998).

$$f_{R_n^{'}}(r_n^{'}) = \frac{(-A^{'}(r_n^{'}))B^n(\theta)}{\Gamma(n)}(A(r_n^{'}))^{n-1}\exp\{-A(r_n^{'})B(\theta)\}$$

therefore

$$\hat{\theta}_{MLE,(R_1^{'},...,R_n^{'})} \overset{Distribution}{=} B^{-1}(\frac{n}{T}) \quad (5)$$

in which

$$T = A(R_n^{'}) \overset{Distribution}{=} Gamma(n, \frac{1}{B(\theta)}) .$$

Considering the above descriptions as well as the relations (4) and (5), the following theorem is readily concluded.

**Theorem 1.** Considering a random sample of size n from the exponential family of the second type and also considering n values from the lower records for this sub family, the below relations are always true.

a) $\hat{\theta}_{MLE,(X_1,...,X_n)} \overset{Distribution}{=} \hat{\theta}_{MLE,(R_1^{'},...,R_n^{'})}$

b) $\gamma(\hat{\theta}_{MLE,(X_1,...,X_n)}) \overset{Distribution}{=} \gamma(\hat{\theta}_{MLE,(R_1^{'},...,R_n^{'})})$ in which $\gamma$ is a real function.

c) $MSE(\gamma(\hat{\theta}_{MLE,(X_1,...,X_n)})) = MSE(\gamma(\hat{\theta}_{MLE,(R_1^{'},...,R_n^{'})}))$

d) $$\alpha^{'}(n) = MSE(\gamma(\hat{\theta}_{MLE,(X_1,...,X_n)})) = MSE(\gamma(\hat{\theta}_{MLE,(R_1^{'},...,R_n^{'})})) =$$
$$E[\gamma\{B^{-1}(\frac{n}{T})\}^2] - 2\gamma(\theta)E[\gamma\{B^{-1}(\frac{n}{T})\}] + \gamma^2(\theta).$$

in which $T \overset{Distribution}{=} Gamma(n, \frac{1}{B(\theta)})$.

## Proof:

Proofs for all parts are obvious considering the relations (4) and (5).

**Note**: From part *d* it is simply found out if a random sample of size *n* from exponential family of the second type includes *m* values of the lower records ( $m \le n$ ), then

$$MSE(\hat{\theta}_{MLE,(R_1^{'},...,R_m^{'})}) = MSE(\hat{\theta}_{MLE,(X_1,...,X_m)}) = \alpha^{'}(m)$$

Simply speaking, the Mean Square Error of the estimations based on m values of lower records is equal with the MSE's of the estimations based on a random sample of size m.

In fact, this is a practically useful note because by considering this note it can be shown that in some cases estimations abased on lower records have less error in comparison with the ones based on random sample. The following examples clarify this

## Example 1.

Let $A(x) = -\log(x), B(\theta) = \theta$, and $\gamma(t) = t$. It is wanted to estimate parameter $B(\theta) = \theta$ in this case $\alpha^{'}(n)$ is obtained as bellow

$$\alpha^{'}(n) = MSE(\gamma(\hat{\theta}_{MLE,(X_1,...,X_n)})) = MSE(\gamma(\hat{\theta}_{MLE,(R_1^{'},...,R_n^{'})})) =$$

$$E[\gamma^2\{B^{-1}(\frac{n}{T})\}] - 2\gamma(\theta)E[\gamma\{B^{-1}(\frac{n}{T})\}] + \gamma^2(\theta)$$

$$= E[\{B^{-1}(\frac{n}{T})\}^2] - 2\theta E[B^{-1}(\frac{n}{T})] + \theta^2 = E[(\frac{n}{T})^2] - 2\theta E[\frac{n}{T}] + \theta^2 =$$

$$\int_0^{+\infty} \frac{n^2}{t^2} \frac{t^{n-1}\theta^n \exp(-\theta t)}{\Gamma(n)} dt - 2\theta \int_0^{+\infty} \frac{n}{t} \frac{t^{n-1}\theta^n \exp(-\theta t)}{\Gamma(n)} dt + \theta^2$$

After some mathematical operations $\alpha^{'}(n)$ is easily obtained as:

$$\alpha^{'}(n) = (\frac{n^2}{(n-2)(n-1)} - \frac{2n}{n-2} + 1)\theta^2 .$$

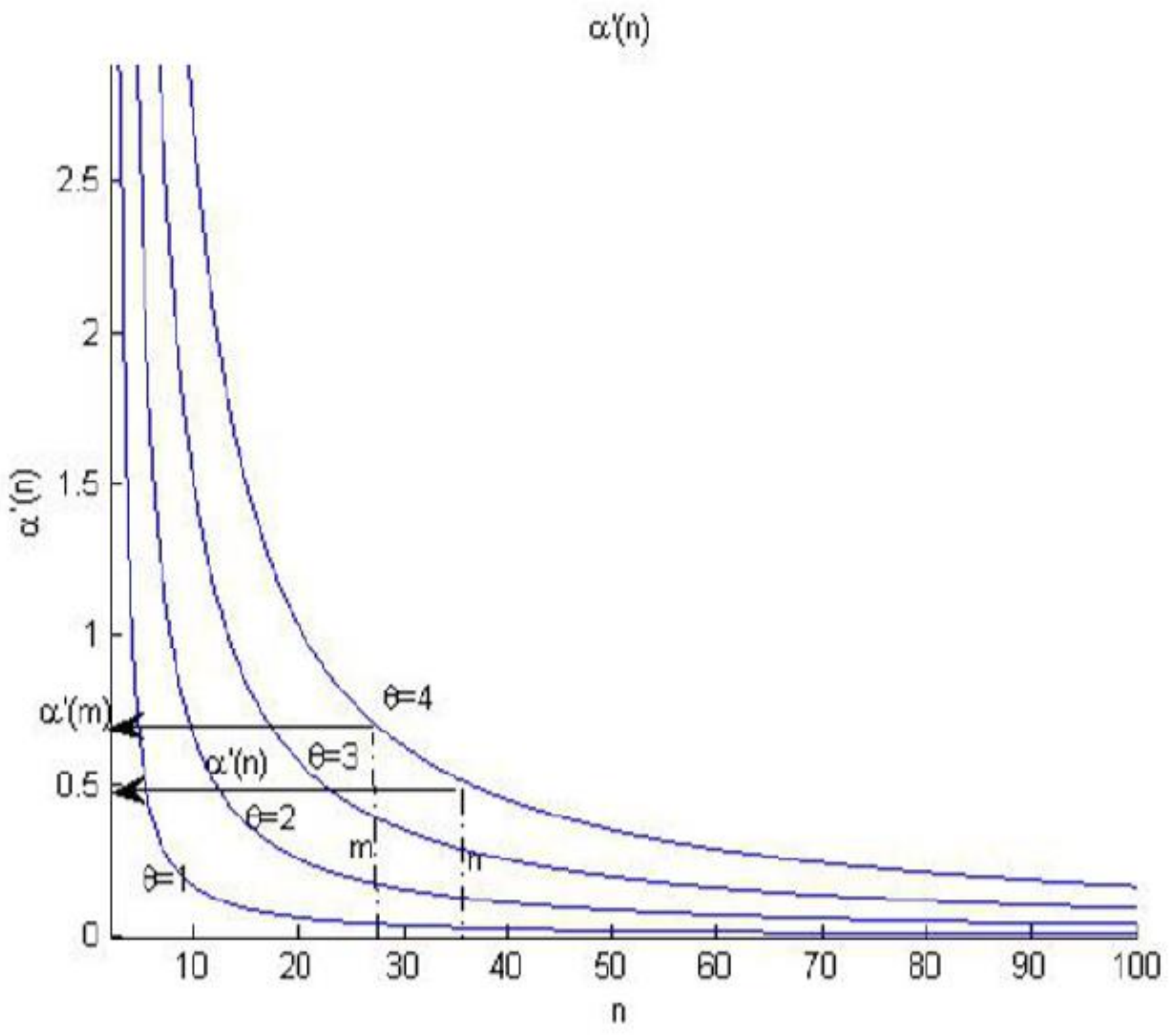

***Figure 1-2 Values of MSE for different values of n.***

Considering the related figure of the above function, it is seen in each sample of size n with m values of lower records $(m < n)$, the estimation based on sample has less value of error in comparison with estimation based on lower record values. It may instill the idea of uselessness of the estimation based on lower record values in the mind, but next instance shows it is not generally true. In some cases, the estimation based on lower record values has less value of error (MSE).

## Example 2.

If in the previous example $\gamma(t) = e^{t}$ , in this case

$$\alpha'(n) = E[\gamma^{2}\{B^{1}(\frac{n}{T})\}] - 2\gamma(\theta)E[\gamma\{B^{-1}(\frac{n}{T})\}] + \gamma^{2}(\theta) =$$

$$E[\exp(2\frac{n}{T})] - 2\exp(\theta)E[\exp(\frac{n}{T})] + \exp(2\theta)$$

$$\alpha'(n) = \int_{0}^{+\infty} \frac{e^{\frac{2n}{t}} t^{n-1}\theta^{n}e^{-\theta t}}{\Gamma(n)} dt - 2e^{\theta} \int_{0}^{+\infty} \frac{e^{\frac{n}{t}} t^{n-1}\theta^{n}e^{-\theta t}}{\Gamma(n)} dt + e^{2\theta}$$

Considering the expansions $e^{\frac{2n}{t}} = \sum_{i=0}^{+\infty} \frac{(2n)^i}{t^i i!}$, and the above relation is obtained as fallow

$$\alpha'(n) = \int_0^{+\infty} \frac{(\sum_{i=0}^{+\infty} \frac{(2n)^i}{t^i i!}) t^{n-1}\theta^n e^{-\theta t}}{\Gamma(n)} dt - 2e^{\theta} \int_0^{+\infty} \frac{(\sum_{i=0}^{+\infty} \frac{n^i}{t^i i!}) t^{n-1}\theta^n e^{-\theta t}}{\Gamma(n)} dt + e^{2\theta} .$$

The value of $\frac{(2n)^i}{t^i i!}$ is greater than zero, that's why it is possible to interchange the integration and summation. Therefore

$$\alpha'(n) = \sum_{i=0}^{+\infty} \frac{(2n\theta)^i \Gamma(n-i)}{\Gamma(i+1)\Gamma(n)} \int_0^{+\infty} \frac{t^{n-i-1}\theta^{n-i} e^{-\theta t}}{\Gamma(n-i)} dt - 2e^{\theta} \sum_{i=0}^{+\infty} \frac{(n\theta)^i \Gamma(n-i)}{\Gamma(i+1)\Gamma(n)} \times$$
$$\int_0^{+\infty} \frac{t^{n-i-1}\theta^{n-i} e^{-\theta t}}{\Gamma(n-i)} dt + e^{2\theta}$$
$$= \sum_{i=0}^{+\infty} \frac{(2n\theta)^i \Gamma(n-i)}{\Gamma(i+1)\Gamma(n)} - 2e^{\theta} \sum_{i=0}^{+\infty} \frac{(n\theta)^i \Gamma(n-i)}{\Gamma(i+1)\Gamma(n)} + e^{2\theta} = \sum_{i=0}^{+\infty} \frac{(n\theta)^i \Gamma(n-i)}{\Gamma(i+1)\Gamma(n)} (2^i - 2e^{\theta}) + e^{2\theta}.$$
$$\Rightarrow \alpha'(n) = \sum_{i=0}^{+\infty} \frac{(n\theta)^i \Gamma(n-i)}{\Gamma(i+1)\Gamma(n)} (2^i - 2e^{\theta}) + e^{2\theta}$$

Considering the graph of obtained $\alpha'(n)$, it can be said if in a sample of size 5, 6, 7, or 8 there is just 2 value of records the MSE of ML estimation based on lower record values has less value than ML estimation based on a random sample. Also, it is clear for samples with large size the MSE of Estimators based on lower record values and based on a random sample are approximately the same. For example, for a sample of size 500 includes 15 value of records, we can say

MSE based on lower records $\overset{approximately}{=}$ MSE based on a random of size 500.

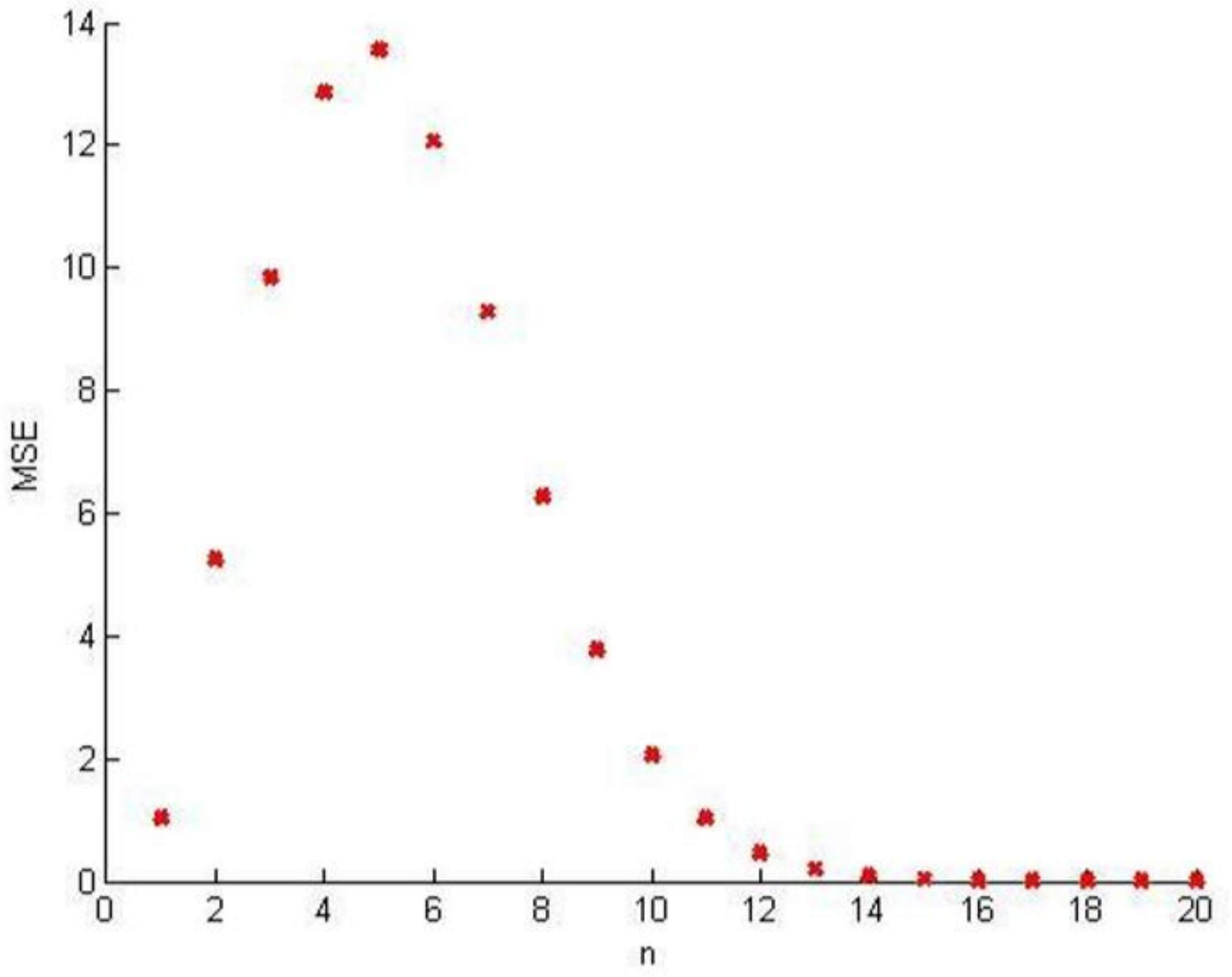

***Figure 2-2. Values of MSE considering different values of n in estimation of the parameter*** $e^{\theta}$ .

## 3- Maximum Likelihood Estimations for Cumulative Distribution Function (CDF) and Probability Density Function (PDF) from Exponential Family of the Second Type Based on Lower Record Values

Regarding stability property of the maximum likelihood estimations and equation (5), hereunder results are readily obtained for PDF and CDF of the type II exponential family.

$$\hat{\theta}_{MLE,(X_1,\dots,X_n)} = B^{-1}\left(\frac{n}{\sum_{i=1}^{n} A\left(X_i\right)}\right) \Rightarrow$$

$$\hat{f}_{MLE,(X_1,\dots,X_n)} = \frac{(-nA^{'}(x))}{\sum_{i=1}^{n} A\left(X_i\right)} \exp\left(-\frac{nA(x)}{\sum_{i=1}^{n} A(X_i)}\right) \overset{Distribution}{=} \frac{(-nA^{'}(x))}{T} \exp\left(-\frac{nA(x)}{T}\right),$$

and

$$\hat{F}_{MLE,(X_1,\dots,X_n)} = \exp\left(-\frac{nA(x)}{\sum_{i=1}^{n} A(X_i)}\right) \overset{Distribution}{=} \exp\left(-\frac{nA(x)}{T}\right).$$

In which $T \overset{Distribution}{=} Gamma(n, \frac{1}{B(\theta)})$ .

**Theorem 2.** If a random sample of size n $(X_1, \dots, X_n)$, includes m values from the lower record values $R_i^{'}, i = 1, \dots, m$, then the maximum likelihood estimations based on lower record values for PDF $f(x;\theta)$ and CDF $F(x;\theta)$ are biased estimators and

$$i - E[\hat{f}_{MLE,(R_1^{'},\dots,R_m^{'})}] = \sum_{i=0}^{+\infty} [\frac{\Gamma(m-i-1)}{\Gamma(m)\Gamma(i+1)} (-1)^i \{mB(\theta)\}^{i+1} (-A^{'}(x))]\{A(x)\}^i$$

$$ii - E[\hat{F}_{MLE,(R_1^{'},\dots,R_m^{'})}] = \sum_{i=0}^{+\infty} [\frac{\Gamma(m-i)}{\Gamma(m)\Gamma(i+1)} \{-mB(\theta)A(x)\}^i ]$$

Here, the second part is proven. The other part is proven similarly.

$$E[\overset{\wedge}{F}_{MLE,(R_1^{'},\dots,R_n^{'})}] = E[\exp(-\frac{mA(x)}{A(R_m^{'})})] = \int_0^{+\infty} \frac{e^{-\frac{mA(x)}{t}} t^{m-1} B^m(\theta) e^{-B(\theta)t}}{\Gamma(m)} dt \qquad (6)$$

Considering the $e^{-\frac{mA(x)}{t}} = \sum_{i=0}^{+\infty} \frac{(-mA(x))^i}{t^i i!}$ and (6)

$$E[\hat{F}_{MLE,(R_1',\dots,R_n')}] = \int_0^{+\infty} \frac{e^{-\frac{mA(x)}{t}} t^{m-1} B^m(\theta) e^{-B(\theta)t}}{\Gamma(m)} dt = \int_0^{+\infty} \frac{\sum_{i=0}^{+\infty} \frac{(mA(x))^i}{t^i i!} t^{m-1} B^m(\theta) e^{-B(\theta)t}}{\Gamma(m)} dt =$$

$$\int_0^{+\infty} \frac{\sum_{even-indices}^{+\infty} \frac{(mA(x))^i}{t^i i!} t^{m-1} B^m(\theta) e^{-B(\theta)t}}{\Gamma(m)} dt - \int_0^{+\infty} \frac{\sum_{odd-indices}^{+\infty} \frac{(mA(x))^i}{t^i i!} t^{m-1} B^m(\theta) e^{-B(\theta)t}}{\Gamma(m)} dt.$$

Considering $A = \int_0^{+\infty} \frac{\sum_{even-indices}^{+\infty} \frac{(mA(x))^i}{t^i i!} t^{m-1} B^m(\theta) e^{-B(\theta)t}}{\Gamma(m)} dt$ and

$$B = \int_0^{+\infty} \frac{\sum_{odd-indices}^{+\infty} \frac{(mA(x))^i}{t^i i!} t^{m-1} B^m(\theta) e^{-B(\theta)t}}{\Gamma(m)} dt.$$

In above expression $\sum_{even-indices}^{+\infty} \frac{(mA(x))^i}{t^i i!}$ and $\sum_{odd-numbers}^{+\infty} \frac{(mA(x))^i}{t^i i!}$ in A and B are positive quantities and that's why it is possible to interchange the integration and summation.

$$E[\hat{F}_{MLE,(R_1',\dots,R_n')}] = \sum_{even-indices}^{+\infty} \frac{\Gamma(m-i)(mA(x)B(\theta))^i}{\Gamma(i+1)\Gamma(m)} \int_0^{+\infty} \frac{t^{m-i-1} B^{m-i}(\theta) e^{-B(\theta)t}}{\Gamma(m-i)} dt -$$

$$\sum_{odd-indices}^{+\infty} \frac{\Gamma(m-i)(mA(x)B(\theta))^i}{\Gamma(i+1)\Gamma(m)} \int_0^{+\infty} \frac{t^{m-i-1} B^{m-i}(\theta) e^{-B(\theta)t}}{\Gamma(m-i)} dt = \sum_{even-indices}^{+\infty} \frac{\Gamma(m-i)(mA(x)B(\theta))^i}{\Gamma(i+1)\Gamma(m)} -$$

$$\sum_{odd-indices}^{+\infty} \frac{\Gamma(m-i)(mA(x)B(\theta))^i}{\Gamma(i+1)\Gamma(m)} = \sum_{i=0}^{+\infty} \frac{\Gamma(m-i)(mA(x)B(\theta))^i}{\Gamma(i+1)\Gamma(m)}.$$

And finally,

$$E[\hat{F}_{MLE,(R_1',R_2',\dots,R_m')}] = \sum_{i=0}^{+\infty} \frac{\Gamma(m-i)}{\Gamma(i+1)\Gamma(m)} (-mA(x)B(\theta))^i.$$

**Theorem 3**. Assuming that a sample of size n from distribution of the exponential family of the second type includes m values of the lower records ($R_1', R_2', \dots, R_m'$), the hereunder results are always true.

$$i-MSE[\hat{f}_{MLE,(R_1',R_2',\dots R_m')}]=\sum_{i=0}^{m-2}(mB(\theta)A'(x))^2[\frac{\Gamma(m-i-2)}{\Gamma(m)\Gamma(i+1)}\{-2mB(\theta)A(x)\}^i]+$$

$$\frac{2}{m}me^{-A(x)B(\theta)}\sum_{i=0}^{m-1}(mA'(x)B(\theta))^2[\frac{\Gamma(m-i-1)}{\Gamma(m)\Gamma(i+1)}\{-mB(\theta)A(x)\}^i]+(A'(x)B(\theta))^2e^{-2A(x)B(\theta)}$$

$$ii-MSE[\hat{F}_{MLE,(R_1',R_2',\dots R_m')}]=\sum_{i=0}^{m-1}\frac{(mA(x)B(\theta)\Gamma(m-i))^i}{\Gamma(i+1)\Gamma(m)}[2^i-2e^{-A(x)B(\theta)}]+e^{-2A(x)B(\theta)}.$$

The second part is proven,

$$MSE[\hat{F}_{MLE,(R_1',R_2',\dots R_m')}]=E[\hat{F}_{MLE,(R_1',R_2',\dots R_m')}-E(\hat{F}_{MLE,(R_1',R_2',\dots R_m')})]^2=$$

$$E(\{\hat{F}_{MLE,(R_1',R_2',\dots R_m')}\}^2)-2F(x;\theta)E(\hat{F}_{MLE,(R_1',R_2',\dots R_m')})+F^2(x;\theta)\Rightarrow$$

$$MSE[\hat{F}_{MLE,(R_1',R_2',\dots R_m')}]=\int_0^{+\infty}\frac{e^{\frac{-2mA(x)}{t}}t^{m-1}B^m(\theta)e^{-B(\theta)t}}{\Gamma(m)}dt-2e^{-A(x)B(\theta)}\int_0^{+\infty}\frac{e^{\frac{-mA(x)}{t}}t^{m-1}B^m(\theta)e^{-B(\theta)t}}{\Gamma(m)}$$

Considering the expansion $e^{-\frac{mA(x)}{t}}=\sum_{i=0}^{+\infty}\frac{(-mA(x))^i}{t^i i!}$,

$$MSE[\hat{F}_{MLE,(R_1',R_2',\dots R_m')}]=\sum_{even-indices}^{+\infty}\frac{(2mA(x))^i}{\Gamma(i+1)\Gamma(m)}\int_0^{+\infty}t^{m-i-1}B^m(\theta)e^{-B(\theta)t}dt-$$

$$\sum_{odd-indices}^{+\infty}\frac{(2mA(x))^i}{\Gamma(i+1)\Gamma(m)}\int_0^{+\infty}t^{m-i-1}B^m(\theta)e^{-B(\theta)t}dt$$

$$2e^{-A(x)B(\theta)}(\sum_{even-indices}^{+\infty}\frac{(mA(x))^i}{\Gamma(i+1)\Gamma(m)}\int_0^{+\infty}t^{m-i-1}B^m(\theta)e^{-B(\theta)t}dt-$$

$$\sum_{odd-indices}^{+\infty}\frac{(mA(x))^i}{\Gamma(i+1)\Gamma(m)}\int_0^{+\infty}t^{m-i-1}B^m(\theta)e^{-B(\theta)t}dt)+e^{-2A(x)B(\theta)}=$$

$$=\sum_{i=0}^{m-1}\frac{(2mA(x))^i\Gamma(m-i)}{\Gamma(i+1)\Gamma(m)}(2^i-2e^{-A(x)B(\theta)})+e^{-2A(x)B(\theta)}$$

In above expression $\sum_{even-indices}^{+\infty}\frac{(2mA(x))^i}{t^i i!}$, $\sum_{odd-indices}^{+\infty}\frac{(2mA(x))^i}{t^i i!}$, $\sum_{even-indices}^{+\infty}\frac{(mA(x))^i}{t^i i!}$, and $\sum_{odd-indices}^{+\infty}\frac{(mA(x))^i}{t^i i!}$ are all positive quantities and that's why it is possible to interchange the integration and summation.

$$MSE[\hat{F}_{MLE,(R_1^{'},R_2^{'},...R_m^{'})}]=\sum_{even-indices}^{+\infty}\frac{(2mA(x))^i}{\Gamma(i+1)\Gamma(m)}\int_0^{+\infty}t^{m-i-1}B^m(\theta)e^{-B(\theta)t}dt-$$

$$\sum_{odd-indices}^{+\infty}\frac{(2mA(x))^i}{\Gamma(i+1)\Gamma(m)}\int_0^{+\infty}t^{m-i-1}B^m(\theta)e^{-B(\theta)t}dt$$

$$2e^{-A(x)B(\theta)}(\sum_{even-indices}^{+\infty}\frac{(mA(x))^i}{\Gamma(i+1)\Gamma(m)}\int_0^{+\infty}t^{m-i-1}B^m(\theta)e^{-B(\theta)t}dt-$$

$$\sum_{odd-indices}^{+\infty}\frac{(mA(x))^i}{\Gamma(i+1)\Gamma(m)}\int_0^{+\infty}t^{m-i-1}B^m(\theta)e^{-B(\theta)t}dt)+e^{-2A(x)B(\theta)}=$$

$$=\sum_{i=0}^{m-1}\frac{(2mA(x))^i\Gamma(m-i)}{\Gamma(i+1)\Gamma(m)}(2^i-2e^{-A(x)B(\theta)})+e^{-2A(x)B(\theta)}$$

**Theorem 4.** Assuming that a sample of size n from the exponential family of the second type $(X_1,...,X_n)$ includes m values of the lower records $(R_1^{'},...,R_m^{'})$, the estimators $\hat{f}_{MLE,(R_1^{'},...,R_m^{'})}$ and $\hat{F}_{MLE,(R_1^{'},...,R_m^{'})}$ are always asymptotically unbiased estimators:

$$i-\lim_{m\to\infty}E[\hat{f}_{MLE,(R_1^{'},...,R_m^{'})}]=-A(x)B(\theta)\exp(-A(x)B(\theta))=f(x;\theta)$$

$$ii-\lim_{m\to\theta}E[\hat{F}_{MLE,(R_1^{'},...,R_m^{'})}]=\exp(-A(x)B(\theta))=F(x;\theta)$$

The following lemma is needed in order to prove the above theorem.

**Lemma 1.** $\lim_{n\to\infty}\frac{\Gamma(n-i-1)n^{i+1}}{\Gamma(n)}=1.$

Because:

$$\lim_{n\to\infty}\frac{\Gamma(n-i-1)n^{i+1}}{\Gamma(n)}=\lim_{n\to\infty}\frac{n^{i}}{\prod_{k=1}^{i+1}\Gamma(n-k)}=\lim_{n\to\infty}\frac{n^{i}}{n^{i}}=1.$$

It must be proven that:

$$\forall \in>0,\exists N,s.t,\forall n\geq N:\frac{\Gamma(n-i-1)n^{i+1}}{\Gamma(n)}-1<\in.$$

Now, assuming that $\in$ is arbitrary, N is chosen as $N=\frac{\sqrt[i+1]{\in+1}-i\sqrt[i+1]{\in+1}}{\sqrt[i+1]{\in+1}-1}$. Now

$$\frac{\Gamma(n-i-1)n^{i+1}}{\Gamma(n)}=\frac{(n-i-2)!n^{i+1}}{(n-1)!}=\frac{n^{i+1}}{(n-1)(n-2)...(n-i-1)}$$
$$=\underbrace{\frac{n}{n-1}\times...\times\frac{n}{n-j-1}}_{(i+1)times}\leq\underbrace{\frac{n}{n-j-1}\times...\times\frac{n}{n-j-1}}_{(i+1)times}$$

On the other hand, it can be simply shown that

$$\forall N\geq n, and, i\in\{1,2,3,...,N-2\}:\frac{n}{n-i-1}\leq\frac{N}{N-i-1},$$

Thus

$$\underbrace{\frac{n}{n-i-1}\times...\times\frac{n}{n-i-1}}_{(i+1)times}<\underbrace{\frac{N}{N-i-1}\times...\times\frac{N}{N-i-1}}_{(i+1)times},$$

and consequently

$$\forall N\geq n:\frac{\Gamma(n-i-1)n^{i+1}}{\Gamma(n)}<\underbrace{\frac{N}{N-i-1}\times...\times\frac{N}{N-i-1}}_{(i+1)times}$$

Finally, by substituting $N=\frac{\sqrt[i+1]{\in+1}-i\sqrt[i+1]{\in+1}}{\sqrt[i+1]{\in+1}-1}$ in the last relation, the following result is reached:

$$\frac{\Gamma(n-i-1)n^{i+1}}{\Gamma(n)}<\in+1.$$

## Proof:

Only the first part is proven; the other parts have similar proofs. First notice that

$$\lim_{m\to\infty}[\hat{f}_{MLE,(R_1',...,R_m')}] = \sum_{i=0}^{\infty}\lim_{m\to\infty}\frac{\Gamma(m-i-1)}{\Gamma(m)\Gamma(i+1)}(-1)^i m^{i+1}\{B(\theta)\}^{i+1}A'(x)\{A(x)\}^i =$$

$$\sum_{i=0}^{\infty}\lim_{m\to\infty}\left(\frac{\Gamma(m-i-1)m^{i+1}}{\Gamma(m)}\frac{\{-B(\theta)\}^i B(\theta)A'(x)\{A(x)\}^i}{\Gamma(i+1)}\right.$$

Considering lemma 1

$$\lim_{m\to\infty}[\hat{f}_{MLE,(R_1',...,R_m')}] = \sum_{i=0}^{\infty}\frac{\{-A(x)B(\theta)\}^i}{i!}B(\theta)A'(x) = B(\theta)A'(x)\exp\{-A(x)B(\theta)\}.$$

**Theorem 5.** In the exponential family of the second type, if the sample $(X_1,...,X_n)$ has m values from the lower records $(R_i')_{i=1,...,m}$, then the maximum likelihood estimations for the distribution function and the probability density function based on random sample and on lower records are consistent estimators for the probability density function and for the distribution function. Then:

$$i-\hat{f}_{MLE,(R_1',...,R_m')}\xrightarrow{P} f_\theta(x)$$

$$ii-\hat{F}_{MLE,(R_1',...,R_m')}\xrightarrow{P} F_\theta(x)$$

**Proof.** Here it is just proven the first part; the second part has a similar proof. It is considered the following abbreviation just for easiness.

$$\hat{f}_{MLE,(R_1',...,R_m')} = \hat{f}_{MLE,R'}$$

Based on the Markov's theorem

$$P\{\varepsilon \le \left|\hat{f}_{MLE,R'} - E[\hat{f}_{MLE,R'}]\right|\} \le \frac{\operatorname{var}(\hat{f}_{MLE,R'})}{\varepsilon}.$$

On the other hand, it is concluded by some mathematical calculations that

$$E[\{\hat{f}_{MLE,R'}\}^2] = \sum_{i=0}^{+\infty}\frac{\Gamma(m-i-2)}{\Gamma(i+1)\Gamma(m)}\{A'(x)\}^2\{B(\theta)\}^2(-2)^i\{mB(\theta)\}^i\{A(x)\}^i,$$

and consequently:

$$Var[\hat{f}_{MLE,R'}] = \sum_{i=0}^{+\infty} \frac{\Gamma(m-i-2)}{\Gamma(i+1)\Gamma(m)} \{A'(x)\}^2 \{B(\theta)\}^2 (-2)^i \{mB(\theta)\}^i \{A(x)\}^i - \{E[\hat{f}_{MLE,R'}]\}^2.$$

Thus

$$\lim_{m \to +\infty} Var[\hat{f}_{MLE,R'}] = \sum_{i=0}^{+\infty} Lim[\frac{m^i \Gamma(m-i-2)}{\Gamma(m)}] \frac{\{-2B(\theta)\}^i \{A(x)\}^i}{\Gamma(i+1)} \{A'(x)\}^2 \{B(\theta)\}^2 - \{f_\theta(x)\}^2 =$$
$$\{A'(x)\}^2 \{B(\theta)\}^2 \exp\{-2A(x)B(\theta)\} = 0$$

Therefore

$$\lim_{m \to +\infty} P\{\left|\hat{f}_{MLE,R'} - E[\hat{f}_{MLE,R'}]\right| > \varepsilon\} \le 0 \text{ so } \hat{f}_{MLE,R'} \xrightarrow{P} E[\hat{f}_{MLE,R'}] \text{ or } \hat{f}_{MLE,R'} - E[\hat{f}_{MLE,R'}] \xrightarrow{P} 0.$$

On the other hand, it is reminded that if $\lim_{m \to +\infty} a_m = a$ and $Y \xrightarrow{P} b$ then (Billingsley (1995))

$$Y + a_m \xrightarrow{P} b + a \qquad (7)$$

Considering (7) and respecting the first part of theorem 4 ( $\lim_{m \to +\infty} E[\hat{f}_{MLE,R'}] = f_\theta(x)$), it is obtained that

$$\hat{f}_{MLE,R'} - E[\hat{f}_{MLE,R'}] + E[\hat{f}_{MLE,R'}] - f \xrightarrow{P} 0,$$

or equally $\hat{f}_{MLE,(R'_1,...,R'_m)} \xrightarrow{P} f_\theta(x)$.

**Theorem 6.** If the sample $(X_1,...,X_n)$ has m values from the lower records ( $R'_i$ $i=1,...,m$ ), then the maximum likelihood estimations based on lower records are consistent estimations for $\theta$.

$$\hat{\theta}_{MLE,(R'_1,...,R'_m)} \xrightarrow{P} \theta$$

## Proof:

With respect to the continuity of functions $B(t), -Log(1-t)$, and $F(t) = 1 - \exp\{-A(t)B(\theta)\}$, also the continuity theorem of composite functions, it is concluded that the function $D(t) = B(\frac{-Log(1-F(t))}{A(x)})$ is always continuous on $\mathbb{R}$. On the other hand, it is known from the second part of theorem 5 that $\hat{F}_{MLE,(R_1^{'},\ldots,R_m^{'})} \overset{P}{\rightarrow} F_{\theta}(x)$, or equivalently:

$$A(x)B(\hat{\theta}_{MLE,(R_1^{'},\ldots R_m^{'})})\exp\{-A(x)B(\hat{\theta}_{MLE,(R_1^{'},\ldots R_m^{'})})\} \overset{P}{\rightarrow} F_{\theta}(x).$$

Now

$$D(\hat{F}_{MLE,(R_1^{'},\ldots R_m^{'})}) \overset{P}{\rightarrow} D(F_{\theta}(x)),$$

which is equivalent to

$$\hat{\theta}_{MLE,(R_1^{'},\ldots R_m^{'})} \overset{P}{\rightarrow} \theta.$$

**4. Conclusions**

In the beginning, a general subfamily from the exponential family was considered, it was tried to obtain the ML estimations based on a random sample and based on lower record values for unknown parameter of this family; then, relations between these estimations were shown theoretically. It was engaged via some theorems to description of this problem that the record estimations can be acceptable in general state. Additionally, the ML estimations based on a random sample and based on lower record values were calculated for the CDF and PDF of this sub-family, and their corresponding errors were described and studied by using some theorems; in addition, it was proven that at asymptotical states, these estimations are unbiased. Subsequently, several theorems associated with convergence in probability were expressed for the resulted estimators.